\documentclass[12pt, twoside, leqno]{article}
\usepackage{latexsym}
\usepackage{amsmath}
\usepackage{amsfonts}
\usepackage{amssymb}
\usepackage{amsthm}
\usepackage{color}
\numberwithin{equation}{section}

\input xy
\xyoption{all} \hfuzz=1.5pt \tolerance=500
\theoremstyle{plain}

\theoremstyle{definition}

\newcommand{\cd}{{\cdot}}

\newcommand{\ot}{{\otimes}}

\newcommand{\va}{\varphi}
\newcommand{\sq}{\square}

\newcommand{\ii}{\infty}

\newcommand{\nl}{\nu\in\Lambda}
\newcommand{\co}{{\mathbb C}}

     \newcommand{\cc}{{\cal CC}}

     \newcommand{\lm}{{\lambda}}

\def\x{\relax\ifmmode {\mbox{*}}\else*\fi}

\newcommand{\ed}{\end{document}}

\begin{document}
\pagestyle{myheadings} \markboth{\centerline{\small{\sc
            A.~Ya.~Helemskii}}}
         {\centerline{\small
         {\sc Projective and free matricially normed spaces}}} 
%Proto-quantum spaces and their tensor products}}}
 \title{\bf Projective and free matricially normed spaces}
 \footnotetext{This research was supported by the Russian
Foundation for Basic Research (grant No. 15-01-08392).}

\author{\bf A.~Ya.~Helemskii}

\date{}

\maketitle

\begin{abstract}
We study metrically projective and metrically free matricially normed spaces.
 We describe these spaces in terms of a special space $\widehat M_n$, the space of $n\times n$
 matrices, endowed with a special matrix--norm. We show that metrically free matricially normed spaces
 are matricial $\ell_1$--sums of some distinguished families of matricially normed spaces $\widehat M_n$,
 whereas metrically projective matricially normed spaces are complete direct summands of matricial 
 $\ell_1$--sums of arbitrary families of the spaces $\widehat M_n$. At the end we specify the 
 underlying normed space of $\widehat M_n$ and show that the spaces $\widehat M_n;n>1$ do not belong to 
 any of the classes ${L}^p; p\in[1,\ii]$, introduced by Effros and Ruan. However, in a certain sense the 
 behavior  of $\widehat M_n$ resembles that of ${L}^1$--spaces.
 \end{abstract}

\bigskip
\centerline{{\bf 1. Introduction}}

\bigskip
The notion of a matricially normed space was introduced by Effros and Ruan~\cite{er1}. Soon, after 
the discovery of Ruan representation theorem~\cite{ru2}, 
 the most attention was concentrated to the outstanding special class 
of these structures, namely the ${L}^\ii$--spaces, more often now called (abstract) operator 
spaces. 
 (see the textbooks~\cite{efr,paul,pis,blem}). However,
 already in~\cite{er1,ru2} it was demonstrated that  matricially normed spaces are a subject of 
considerable interest also outside the class of operator spaces. In particular, according to these 
papers, one can successfully study the Haagerup tensor product of these spaces. Later it was shown 
in~\cite{colh} that another important tensor product of operator spaces, the projective tensor 
product, 
%of operator spaces, 
also can be studied in the context of general matricially normed spaces. 

In this paper we introduce and study the closely connected notions of a metrically projective and a 
metrically free matricially normed space. In the realm of operator spaces the definition of a 
metrically projective space resembles what was called by Blecher~\cite{ble} (just) projective  
space, but differs from the latter. In the classical context of Banach spaces the metric 
projectivity appeared, under a different name, in the old paper of Graven~\cite{grav}. 

We first characterize metrically free matricially normed spaces. Then, using the well known general 
categorical connection between freeness and projectivity, we characterize metrically projective 
matricially normed spaces. We describe the spaces in both classes in terms of the special space 
$\widehat M_n$, the space of $n\times n$--matrices endowed with a special matrix--norm. We show 
that metrically free spaces are matricial $\ell_1$-sums of some specified families of spaces 
$\widehat M_n$, whereas metrically projective spaces are complete direct summands of matricial 
$\ell_1$-sums of arbitrary families of spaces $\widehat M_n$. 

The latter result, concerning projectivity, resembles Blecher's Theorem 3.10 in~\cite{ble}, and in 
fact we were inspired by that. Note that the mentioned theorem, extended in a straightforward way 
to general matricially normed spaces, also can be deduced from the description of metrically free 
matricially normed spaces, however after some elaboration of our general categorical tools; 
cf.~\cite{he6}. 
%(see the `asymptotic projectivity', discussed in~\cite{???}). 
But we leave this material outside the scope of the present paper.  

The contents of the paper are as follows.

Section 2 contains some preliminary definitions.

In Section 3 we introduce our main matricially normed space $\widehat M_n$. 

In Section 4 we prepare our tools from category theory. We consider the so-called rigged categories 
(well known under many different names), define projective objects in such a category and show that 
the metric projectivity of matricially normed spaces is a particular case of this general 
categorical projectivity. Then we introduce, within the frame--work of a rigged category, the 
notion of a free object.
% and, as a particular case, that of a metrically free matricially normed space. 
We recall several general categorical observations that will be used in later sections, notably the 
characterization of projective objects as retracts of free objects.  

In Section 5 we fix $n\in{\Bbb N}$ and introduce the special rig `$\odot_n$', playing, in a sense, 
the role of a `building brick' for 
%our main rig `$\odot$', the latter 
the rig, responsible for the metric projectivity. We show that the $\odot_n$--free object with a 
one-point base 
%, corresponding to the rig `$\odot_n$', 
is exactly the matricially normed space $\widehat M_n$. The proof heavily relies on properties of a 
distinguished $n^2\times n^2$--matrix ${\Bbb I}$. The latter, in the guise of 
% that, being identified with 
an element of $M_n\ot M_n$, is $\sum_{ij}e_{ji}\ot e_{ij}$, where $e_{ij}$ denotes the elementary 
matrix with 1 on the $ij$--th place. 

In Section 6 we apply the results of the previous sections to obtain our main results; namely, the 
above--mentioned description of metrically free and (as a corollary) of metrically projective 
matricially normed spaces. 

Finally, in Section 7 we obtain some (far from complete) information about the structure of the 
space $\widehat M_n$. First, we find its underlying normed space: it turns out to be the space of 
$n\times n$--matrices with the trace class norm. Thus, it is the same as the underlying space of 
the operator space $T_n$, playing the main role in the description of projective {\it operator} 
spaces in~\cite{ble}. However, as a matricially normed space, $\widehat M_n$ is profoundly 
different from $T_n$. We show that 
% (not surprisingly) 
it does not belong to any  of the classes $L^p; 1<p\le\ii$ of Effros/Ruan, which is not surprising, 
and (what is somehow surprising to the author) to the class $L^1$ as well. Nevertheless, in a 
certain sense the behavior of $\widehat M_n$ resembles that of ${L}^1$--spaces. 

\bigskip
 \centerline{{\bf 2. Initial definitions}}

\bigskip
In what follows, we denote the space of $n\times n$-matrices, as a pure algebraic object, by  
$M_n$, and the same space, endowed by  the operator norm $\|\cd\|_o$ or the trace norm $\|\cd\|_t$, 
by ${\Bbb M}_n$ and ${\Bbb T}_n$, respectively. If $E$ is a normed space, we denote by $B_E$ its 
closed unit ball. The identity operator on a linear space $E$ is 
 denoted by ${\bf 1}_E$.

Let $E$ be a linear space, $M_n(E)$ the space of $n\times n$-matrices with entries from $E$. We  
identify $M_n(E)$ with the tensor product $M_n\ot E$. According to our convenience, we shall use 
either `matrix guise' or `tensor guise' of this space. We denote the $E$--valued diagonal 
block--matrix with matrices $u_1,...,u_k$ on the diagonal by $u_1\oplus...\oplus u_k$.  

\medskip
{\bf Definition 1.} (Effros/Ruan~\cite{efr}) A sequence of norms $\|\cd\|_n$ on $M_n(E): n\in{\Bbb 
N}$ is called a {\it matrix--norm} on $E$, if it satisfies the two following conditions: 

%\medskip
{\bf Axiom 1}. For $u\in M_n(E)$ and $n,m\in{\Bbb N}$ we have $\|u\oplus0\|_{n+m}=\|u\|_n$. Here 0 
is the zero matrix in $M_m(E)$. 

%\medskip
{\bf Axiom 2.} For $u\in M_n(E)$ and $S\in M_n$ we have $\|Su\|_n\le\|S\|_o\|u\|_n$ and 
$\|uS\|_n\le\|u\|_n\|S\|_o$. 

\medskip
A space $E$, endowed with a matrix--norm, is called {\it  matricially normed space}. The normed 
space, identified with $M_1(E)$, is called {\it underlying normed space} of our  matricially normed 
space. 
%(Quantization?) 

Two examples that we shall need are the matricially normed spaces $\co_{\min}$ and $\co_{\max}$. 
The first one is $\co$ with the matrix-norm, arising after the identifying, for every $m\in{\Bbb 
N}$, of $M_m(\co)$ with ${\Bbb M}_m$, whereas the second one is $\co$ with the matrix-norm, arising 
after the identifying of $M_m(\co)$ with ${\Bbb T}_m$. 

Every subspace $F$ of a matricially normed space is, of course, a  matricially normed space with 
respect to induced norms in $M_n(F)\subseteq M_n(E)$ for every $n$. We call it  matricially normed 
subspace of $E$. 

\bigskip
Now let $E$ and $F$ be linear spaces, $\va:E\to F$ a linear operator. The operator $\va_n:M_n(E)\to 
M_n(F):(x_{ij})\mapsto(\va(x_{ij}))$, is called {\it $n$--th amplification} of $\va$. (In the 
`tensor approach' $\va_n$ is, of course, ${\bf 1}_{M_n}\ot\va$). 

If our $E$ and $F$ are matricially normed spaces, then we call $\va$ {\it completely bounded}, if 
$\sup\{\|\va_n\|;n\in{\Bbb N}\}<\ii$. We denote this supremum by $\|\va\|_{cb}$. 
 
If, in the previous context, every amplification is a contractive operator (that is 
$\|\va\|_{cb}\le1$), we say that $\va$ is {\it completely contractive}. (This is the most important 
class of operators in the present paper). The set of completely contractive operators between $E$ 
and $F$ is denoted by $\cc(E,F)$. If every amplification is isometric, strictly coisometric or 
isometric isomorphism, we say that $\va$ is {\it completely isometric, completely strictly 
coisometric or completely isometric isomorphism}, respectively. Here we recall that the operator 
between normed spaces $E$ and $F$ is called strictly coisometric (or exact quotient map), if it 
maps $B_E$ {\it onto} $B_F$. 

The (non--additive) category with 
%$MN$--пространства 
matricially normed spaces as objects and completely contractive operators as morphisms is denoted 
by ${\bf MN_1}$. Evidently, isomorphisms in this category are completely isometric isomorphisms, 
defined above. 

\bigskip
{\bf Definition 2}. A matricially normed space $P$ is called {\it metrically projective}, if, for 
every completely strictly coisometric operator $\tau$ between matricially normed spaces, say $E$ 
and $F$, and every completely contractive operator $\va:P\to F$ there exists a completely 
contractive operator $\psi:P\to E$, making the diagram 
$$
\xymatrix@R-10pt@C+15pt{
& E \ar[d]^{\tau}\\
P \ar[ur]^{\psi} \ar[r]^{\va} & F }\eqno(D1)
$$
\noindent commutative. 

\bigskip
\centerline{{\bf 3. Construction of the matricially normed spaces $\widehat M_n$}} 

\bigskip 
From now on and until we state otherwise, we fix some $n\in{\Bbb N}$. Sometimes, if there is no 
danger of confusion, we omit index $n$. 

We begin with pure algebraic preparations. 
     
Suppose we are given a linear space $E$. Let us introduce the operator 
$$
\iota^E:M_n(E)\to{\cal L}(M_n,E),\eqno(1)
$$
where ${\cal L}(\cd,\cd)$ is the symbol of the space of linear operators. It takes an $n\times 
n$--matrix $(x_{ij})$ with entries in $E$ to the operator 
$(\lm_{ij})\mapsto\sum_{ij}\lm_{ji}x_{ij}\in E$, where $(\lm_{ij})\in M_n$. Equivalently, if we use 
the `tensor guise' of $M_n(E)$, then $\iota^E:M_n\ot E\to{\cal L}(M_n,E)$ is well defined by taking 
an elementary tensor $a\ot x$ to the operator $b\mapsto tr(ab)x$. Here and onwards $tr$ denotes the 
trace of a matrix. Obviously, $\iota^E$ is a linear isomorphism. 
 
It is convenient to denote, for a given $v\in M_n(E)$, the operator $\iota^E(v)$ by $\va^v:M_n\to 
E$. 

Now consider all possible couples $(E,v)$, where $E$ is a matricially normed space 
 and $v\in B_{M_n(E)}$. In what follows, we refer them as {\it proper couples.} 

\medskip
{\bf Definition 3}. For every $m\in{\Bbb N}$ and $u\in M_m(M_n)$ we set 
\[
\|u\|_m:=\sup\{\|(\va^v)_m(u)\|_m\},\eqno(2)
\]
 where supremum is taken over all proper couples. 
 
 \medskip
 (We recall that $(\va^v)_m:M_m(M_n)\to M_m(E)$ takes an $m\times m$--matrix $(a^{kl})$ with 
 entries in $M_n$ to the $m\times m$--matrix $(\va^v(a^{kl})$) with entries in $E$.
 
\medskip
{\bf Proposition 1.} {\it The indicated supremum is finite. Moreover, $\|u\|_m$ does not exceed the 
sum of modules of matrix entries after the identification of $M_m(M_n)$ with $M_{mn}$. } 

\smallskip
{\it Proof}. Take $u\in M_m(M_n)$; $u=(a^{kl}), a^{kl}=(\lm^{kl}_{ij})\in M_n$. We must show that 
for every proper couple $(E,v)$ we have $\|(\va^v)_m(u)\|\le\sum_{kl}\sum_{ij}|\lm^{kl}_{ij}|$. 

If $v=(x_{ij});x_{ij}\in E$, then, since $v\in B_{M_n(E)}$, and $E$ is a 
%n $MN$--space
 matricially normed space, we have $\|x_{ij}\|\le1$ for all $i,j$. Hence 
 $\|\va^v(a^{kl})\|\le\sum_{ij}|\lm^{kl}_{ij}|$. On the other hand, using again that $E$ is a
%n $MN$--space
 matricially normed space, we have $\|(va^v)_m(u)\|\le\sum_{kl}\|\va^v(a^{kl})\|$. $\sq$

\medskip
{\bf Theorem 1.} {\it The sequence of functions $u\mapsto\|u\|_m;m\in{\Bbb N}$ is a matrix-norm on 
$M_n$. } 

\smallskip
{\it Proof}. First, we show that the function $u\mapsto\|u\|_m$ is a seminorm on $M_m(M_n)$ and 
then we check the Axioms 1 and 2. All three assertions are proved by a similar argument. Namely, we 
use the respective properties of matricially normed spaces $E$ in proper couples and then the 
definition of $\|u\|_m$ as the relevant supremum. For example, the first estimation in Axiom 2 
follows from the relations 
$$
\|(\va^v)_m(Su)\|=\|S(\va^v)_m(u)\|\le\|S\|_o\|(\va^v)_m(u)\|\le\|S\|_o\|\|u\|_m.
$$

Finally, we prove that our seminorm is actually a norm. Take $u\in M_m(M_n);u= 
% as a block--matrix 
(a^{kl}); a^{kl}=(\lm^{kl}_{ij})\in M_n$. If $u\ne0$, then $\lm^{k'l'}_{j'i'}\ne0$ for a certain 
$k',l'$ and $i',j'$. Choose a proper couple with $v=(x_{ij})\in M_n(E)$ such that $x_{ij}\ne0$ if, 
and only if $i=i'$ and $j=j'$, and also $\|x_{i'j'}\|\le1$. Then we see that 
 $(\va^v)_m$ takes $u$ to the matrix $(y_{kl})\in M_m(E)$ with $y_{k'l'}\ne0$. Therefore
 $\|(\va^v)_m(u)\|>0$, hence $\|u\|>0$. $\sq$

\medskip
We denote the resulting matricially normed space by $\widehat M_n$. 

In particular, it is easy to show that $\widehat M_1$ is just $\co_{\max}$. The structure of 
$\widehat M_n$ for bigger $n$ is not so transparent; we shall see this in our last section. 

\bigskip
\centerline{{\bf 4. Projectivity and freeness in rigged categories}} 

\bigskip
{\bf Definition 4.} Let ${\cal  K}$ be an arbitrary 
%(generally speaking, non-additive) 
category. A {\it rig} of ${\cal  K}$ is a faithful covariant functor $\square:{\cal K}\to{\cal L}$, 
where ${\cal L}$ is another category. A pair $({\cal  K},\sq)$, consisting of a category and its 
rig, is called {\it rigged category}. 

We call a morphism $\tau$ in ${\cal K}$ {\it admissible}, if $\square(\tau)$ is a retraction in 
${\cal L}$. %(Clearly, such a $\tau$ is indeed an epimorphism.) 

\medskip
{\bf Definition 5}. An object $P$ in ${\cal K}$ is called {\it $\square$-projective}, if, for every 
$\sq$--admissible morphism $\tau:Y\to X$ and every morphism $\va:P\to X$, there exists a morphism 
$\psi:P\to Y$, making the diagram (D1) commutative. 

\medskip {\bf Our principal example. } Consider the covariant functor 
\[
\odot:{\bf MN_1}\to{\bf Set},\eqno(3)
\]
 taking a matricially normed space $E$
 %n $MN$--space $E$ 
 to the cartesian product $\textsf{X}_{m=1}^\ii B_{M_m(E)}$. Thus, the elements of the set $\odot(E)$ 
 are sequences 
$(w_1,\dots,w_m,\dots)$ where $w_m\in B_{M_m(E)}$. As to the action of our functor on morphisms, it 
takes a completely contractive operator $\va:E\to F$ to the map 
\[
\odot(\va):\odot(E)\to\odot(F), (w_1,\dots,w_m,\dots)\mapsto(\va(w_1),\dots,\va(w_m),\dots).
\]
It is clear that we have obtained a rigged category. Note the following obvious statement. 

\medskip
{\bf Proposition 2}. {\it (i) A morphism in ${\bf MN_1}$ is $\odot$--admissible 
%contractive operator between matricially normed spaces is a $\odot$--admissible morphism 
if, and only if it is a completely strictly coisometric operator. 

(ii) $\odot$-projective objects in ${\bf MN_1}$ are exactly metrically projective matricially 
normed spaces}. $\sq$ 

\bigskip
Return to general rigged categories. The following concept is well known under different names. 

\medskip
{\bf Definition 6}. Let $M$ be an object in ${\cal L}$. An object ${\bf Fr}(M)$ in ${\cal K}$ is 
called {\it $\square$-free object with the base $M$}, if for every $X\in{\cal K}$, there exists a 
bijection 
$$
{\cal I}_X: {\bf h}_{{\cal L}}(M,\square X)\to{\bf h}_{{\cal K}}({\bf Fr}(M),X),\eqno(4)
$$
between the respective sets of morphisms, natural on $X$. We say that a rigged category {\it admits 
freeness}, if every object in ${\cal L}$ is a base of a free object in ${\cal K}$. 

\medskip
{\bf Remark.} According to~\cite[Chs. III,IV]{mc2}, to say that a rigged category admits freeness 
is equivalent to say that $\sq$ has a left adjoint functor. 

\medskip
 The following observations show the practical use of the freeness. They are actually well known 
and can be extricated, as particular cases or easy corollaries, from some general facts, contained 
in~\cite[Chs. III,IV]{mc2}. 

\medskip
{\bf Proposition 3}. {\it Suppose that our rigged category admits freeness. Then 

(i) every object in ${\cal K}$ is the range of an admissible morphism with a free domain. 

(ii) an object in ${\cal K}$ is is projective if, and only if it is a retract of a free object. In 
particular, all free objects in ${\cal K}$ are projective.} $\sq$ 

\medskip
It was proved in~\cite{he6} that the rig, obtained from $\odot$ by the restriction of ${\bf MN_1}$ 
to its subcategory of {\it operator spaces}, admits freeness, and its free objects are the 
so--called $\oplus_1$--sums
% coproducts {\it in this subcategory} 
of certain families of `building bricks'. The latter are the spaces,
% ${\Bbb T}_n$, 
considered, as 
operator spaces, as dual to the `concrete' operator space ${\Bbb M}_n$. (These spaces were already 
used in~\cite{ble}). A similar result was obtained in~\cite{he7} in more general context of 
operator modules over operator algebras. 

But our aim is to find free objects in the `whole' rigged category $({\bf MN_1},\odot)$. To begin 
with, we shall find free objects in another rigged category that is, speaking informally, a `small 
part' of $({\bf MN_1},\odot)$.
% much simpler than $({\bf MN_1},\odot)$. 

\newpage
 
\bigskip
\centerline{{\bf 5. The rig $\odot_n$ and its free objects with the one-point base}} 

\bigskip
With $n$ still fixed, we consider the covariant functor 
\[
\odot_n:{\bf MN_1}\to{\bf Set},\eqno(5)
\]
taking a matricially normed space $E$
%n $MN$--space $E$ 
to the set $B_{M_n(E)}$ and taking a completely contractive operator to its restriction to the 
respective unit balls. Evidently, we get a rig. 

Now we need some preparation. First, recall the operator $\iota^E$ (see (1)). 

\medskip
{\bf Proposition 4}.{\it Let $E,F$ be linear spaces, $\psi:E\to F$ an operator. Then the diagram 
$$
\xymatrix@C+20pt{M_n(E) \ar[r]^{\iota^E}\ar[d]_{\psi_n}
& {\cal L}(M_n,E) \ar[d]^{h_\psi} \\
M_n(F) \ar[r]^{\iota^F} & {\cal L}(M_n,F), }\eqno(D2)
$$
\noindent where $h_\psi$ takes an operator $\chi$ to the composition $\psi\chi$, is commutative. } 

\smallskip
{\it Proof}. A convenient way to check this is to use the `tensor guise' of $M_n(E)$ and look at 
elementary tensors in $M_n\ot E$. $\sq$ 

\medskip
Recall the matrices $e_{ij}\in M_n$ and ${\Bbb I}\in M_n(M_n)$ from Section 1.
%, which is $\sum_{ij}e_{ji}\ot e_{ij}$ after the identification of $M_n(M_n)$ with $M_n\ot M_n$ (cf. Section 1). 
A routine calculation gives 
%`$ij$-th' place. Set $E:=M_n$ and denote by ${\Bbb I}$ the element $\sum_{ij}e_{ji}\ot e_{ij}$ in the `tensor guise' of $M_n(M_n)$. A routine 

\medskip
{\bf Proposition 5}. {\it The operator $\va^{\Bbb I}$ is just ${\bf 1}_{M_n}$.} $\sq$ 
%, the identity operator on $M_n$. $\sq$ % $\tl\q\tr$ 

\medskip
{\bf Proposition 6}. {\it For every linear space $E$ and $v\in M_n(E)$ we have $(\va^v)_n({\Bbb 
I})=v.$ } 

\smallskip
{\it Proof}. Represent $v$ as $\sum e_{ij}\ot x_{ij}; x_{ij}\in E$. Then we have 
$\va^v(e_{ij})=x_{ji}$, hence $(\va^v)_n({\Bbb 
I})=\sum_{ij}e_{ji}\ot\va^v(e_{ij})=\sum_{ij}e_{ji}\ot x_{ji}=v$. $\sq$ 

\medskip
This, together with (2), implies 

\medskip 
{\bf Proposition 7}. {\it The norm of ${\Bbb I}$ in $M_n(\widehat M_n)$ is 1. } $\sq$ 

\medskip
{\bf Theorem 2.} {\it The matricially normed space $\widehat M_n$ is $\odot_n$--free with a 
one-point set as its base. } 

\smallskip
{\it Proof}. Let $E$ be a matricially normed space, and $\{\star\}$ a one-point set. According to 
(4), we must construct a bijection 
$$
{\cal I}_E^n:{\bf h}_{\bf Set}(\{\star\},\odot_n(E))\to{\bf h}_{\bf MN_1}(\widehat M_n,E),%\eqno(2.2)
$$
\noindent or, equivalently, a bijection ${\cal I}_E^n:B_{M_n(E)}\to\cc(\widehat M_n,E)$, natural on 
$E$. Take $v\in B_{M_n(E)}$ and consider $\va^v$ as an operator between the matricially normed 
spaces $\widehat M_n$ and $E$. Then for every $m$ and $u\in M_m(\widehat M_n)$ we have, by (2),  
that $\|u\|_m\ge\|(\va^v)_m(u)\|$. This means that $\va^v:\widehat M_n\to E$ is completely 
contractive. Thus, the map $\iota^E$ has the well defined restriction to $B_{M_n(E)}$ and 
$\cc(\widehat M_n,E)$. It is this restriction that we choose as ${\cal I}_E^n$. Show that it has 
required properties. 

The commutative diagram (D2), being restricted to the respective unit balls and sets 
$\cc(\cd,\cd)$, demonstrates that our constructed ${\cal I}_E^n$ is natural on $E$. Also ${\cal 
I}_E^n$ is obviously injective. It remains to show that it is surjective. 

Take an arbitrary $\psi:\cc(\widehat M_n,E)$ and consider the diagram  
$$
\xymatrix@C+20pt{B_{M_n(\widehat M_n)} \ar[r]^{{{\cal I}_{\widehat M_n}^n}}\ar[d]_{\psi_n}
& \cc(\widehat M_n,\widehat M_n) \ar[d]^{h_\psi} \\
B_{M_n(E)} \ar[r]^{{\cal I}_E^n} & \cc(\widehat M_n,E), }
$$
the relevant restriction of the diagram (D2) after choosing $\widehat M_n$ as $E$ and $E$ as $F$. 
Now recall that the element ${\Bbb I}$ belongs, by Proposition 7, to $B_{M_n(\widehat M_n)}$, and  
Proposition 5 implies that ${\cal I}_{\widehat M_n}^n({\Bbb I})={\bf 1}_{\widehat M_n}$. But 
$h_\psi({\bf 1}_{\widehat M_n})=\psi$, and our diagram is commutative. Therefore $\psi={\cal 
I}_E^n(\psi_n({\Bbb I}))$, and we are done. $\sq$. 

\bigskip
\centerline{{\bf 6. Characterization of free and projective spaces}} 

\bigskip
To move from the rig $\odot_n$ and its free objects with one-point bases to the `whole' $\odot$ and 
its free objects with arbitrary bases, we need the following well known categorical concept (cf., 
e.g.,~\cite[Ch.2]{faith} or~\cite{mc2}). 

\medskip
Let $X_\nu;\nl$ be a family of objects in an (arbitrary) category ${\cal K}$. We recall that a pair 
$(X,\{i_\nu;\nl\})$, where $X$ is an object, and $i_\nu:X_\nu\to X$ are morphisms in ${\cal K}$, is 
said to be the \emph{coproduct} of this family, if, for every object $Y$ and a family of morphisms 
$\psi_\nu:X_\nu\to Y$ there is a unique morphism $\psi:X\to Y$ such that we have $\psi 
i_\nu=\psi_\nu$ for every $\nl$.  

(We speak about `the' coproduct because it is unique up to a categorical isomorphism, compatible 
with the respective coproduct injections.) 

The mentioned $X$, denoted in a detailed form by $\coprod\{X_\nu;\nl\}$, is referred as the {\it 
coproduct object}, and $i_\nu$'s as the {\it coproduct injections}. 
 The morphism $\psi$ is called the {\it coproduct of the morphisms} $\psi_\nu$ and denoted
 by $\coprod\{\psi_\nu;\nl\}$. 

We say that ${\cal K}$ \emph{admits coproducts}, if every family of its objects has the coproduct. 

Of course, the category ${\bf Set}$ admits coproducts: the coproduct of a family of sets is their 
disjoint union, with obvious coproduct injections. Also it is well known that the category ${\bf 
Nor_1}$ of normed spaces and contractive operators also admits coproducts: the coproduct of a 
family of normed spaces is their (classical) $\ell_1$--sum. 

\medskip
Now suppose we have a family $E_\nu;\nl$ of matricially normed spaces. Consider their algebraic sum 
$E:=\oplus\{E_\nu;\nl\}$ and identify, for every $m\in{\Bbb N}$, the linear spaces $M_m(E)$ and 
$\oplus\{M_m(E_\nu);\nl\}$. Endow every $M_m(E)$ with the norm of the $\ell_1$--sum of normed 
spaces. Then we easily see that we made $E$ a matricially normed space. We call it {\it matricial 
$\ell_1$-sum} of a given family. As an easy corollary of the structure of coproducts in ${\bf 
Nor_1}$, we obtain 
%above-mentioned fact concerning ${\bf Nor_1}$, we obtain 

\medskip
{\bf Proposition 8}. {\it The matricial $\ell_1$-sum of a given family of matricially normed spaces 
is the coproduct of this family in ${\bf MN_1}$ with the natural embeddings $i_\nu:E_\nu\to E$ as 
the coproduct injections. Thus, the category ${\bf MN_1}$ admits coproducts.} $\sq$ 

\medskip
{\bf Remark.} The full subcategory of ${\bf MN_1}$, consisting of {\it operator} spaces, also 
admits coproducts, but the respective construction is necessarily more sophisticated. It was shown 
by Blecher~\cite[Sect. 3]{ble}.
% with the help of another, necessarily more sophisticated construction, has shown that  

\medskip
We turn to $\odot$--free objects that in what follows will be referred as {\it metrically free 
matricially normed spaces}. At first we concentrate on the case of the one--point base. 

\medskip
From now on we `release' $n$. Denote by $\widehat M_\ii$ the matricial $\ell_1$--sum ( = coproduct 
in ${\bf MN_1}$) of the family $\{\widehat M_n; n\in{\Bbb N}\}$. 

\medskip
{\bf Theorem 3.} {\it The metrically free matricially normed  space with a one--point base, say 
$\{\star\}$, does exist, and it is $\widehat M_\ii$}. 

\smallskip
{\it Proof}. Let $E$ be an arbitrary matricially normed space. We must construct a bijection 
$$
{\cal I}_E^\ii:{\bf h}_{\bf Set}(\{\star\},\odot(E))\to{\bf h}_{\bf MN_1}(\widehat M_\ii,E),\eqno(6)
$$ 
\noindent natural on $E$. The first of the indicated sets can be identified with the set of 
sequences ${\bf w}=(w_1,...,w_n,...); w_n\in B_{M_n(E)}$. 

By Theorem 2, for every $n$, after relevant identifications, there exists a bijection ${\cal 
I}_E^n:B_{M_n(E)}\to\cc(\widehat M_n,E)$, taking $w_n$ to the operator $\va^{w_n}$. Thus every 
sequence ${\bf w}$ gives rise to a family of completely contractive operators $\va^{w_n}:\widehat 
M_n\to E$. Denote by $\va^{{\bf w}}:\widehat M_\ii\to E$ the coproduct of these $\va^{w_n}$.    
 
Taking every ${\bf w}$ to $\va^{{\bf w}}$, we obtain, modulo the mentioned identifications, a map 
% map from $\{\star\}$ into $\odot(E)$, identified with the respective ${\bf w}$, to 
%$\va^{{\bf w}}$, we obtain a map 
${\cal I}_E^\ii$ between the sets, indicated in (6). Note that 
% Every ${\cal I}_E^n$ is natural on $E$, and 
 for all completely contractive operators 
$\psi:E\to F$, where $E,F$ are matricially normed spaces, we obviously have 
$\psi(\coprod\{\va^{v_n}\})=\coprod\{\psi\va^{v_n}\}$. Therefore, knowing that every ${\cal I}_E^n; 
n\in{\Bbb N}$ is natural on $E$, we obtain that 
%This easily implies that 
${\cal I}_E^\ii$ is natural on $E$. Finally, since every ${\cal I}_E^n$ is a bijection, ${\cal 
I}_E^\ii$ is also a bijection. $\sq$ 

\medskip
To pass from one--point sets, as bases of free objects, to arbitrary sets, we can use the following 
simple categorical observation. Let $\sq:{\cal K}\to{\cal L}$ be an arbitrary rig. 

\medskip
{\bf Proposition 9.} \emph{Suppose that we are given a family  $F_\nu;\nl$ of free objects with 
bases $M_\nu$. Further, suppose that there exist the coproducts $F:=\coprod\{F_\nu;\nl\}$ and 
$M:=\coprod\{M_\nu;\nl\}$ in ${\cal K}$ and ${\cal L}$, respectively. Then $F$ is a free object 
with the base $M$.} 

\smallskip
{\it Proof}. See, e.g.,~\cite[Prop. 2.13]{he6}. $\sq$ 

\medskip
Since every set $M$ is the coproduct of its one-point subsets, this proposition immediately implies 
%, for the case of  our principal rig $\odot$,

\medskip 
{\bf Theorem 4.} {\it For every set $M$, there exists a metrically free  matricially normed space 
with the base $M$, and it is the matricial $\ell_1$-sum 
%( =coproduct in ${\bf MN_1$) 
of the family of copies of the matricially normed space $\widehat M_\ii$, indexed by points of $M$. 
Thus, the rigged category $({\bf MN_1},\odot)$ admits freeness.} $\sq$  

\medskip
Now we want to pass from free to projective  matricially normed spaces. To make the formulation 
more geometrically transparent, we say that a matricially normed space $F$ is a a complete direct 
summand of a matricially normed space 
% an $MN$--space 
$E$, if $F$ is completely isometrically isomorphic to a  matricially normed subspace $G$ of $E$, 
and there is a completely contractive projection of $E$ onto $G$. We have an obvious

\medskip
{\bf Proposition 10}. {\it A matricially normed space $F$ is a retract in ${\bf MN_1}$ of a 
matricially normed space $E$ if, and only if $F$ is a complete direct summand of $E$}. $\sq$ 

In what follows, we use a simple general-categorical observation, concerning an arbitrary rig. 

\medskip
{\bf Proposition 11.} {\it (i) A retract of a $\sq$--projective object is $\sq$--projective

(ii) the coproduct of a family of $\sq$--projective objects} (if, of course, it does exist) {\it is 
$\sq$--projective.} $\sq$

\medskip
We call a matricially normed space {\it $\widehat M$--composed}, if it is a matricial $\ell_1$--sum 
of some family of spaces such that each of summands is $\widehat M_n$ for some $n\in{\Bbb N}$. 

\medskip
{\bf Theorem 5.} {\it (i) Every matricially normed space is an image of a completely strictly 
coisometric operator with the $\widehat M$--composed space as its domain

(ii) A matricially normed space is metrically projective if, and only if it is a complete direct 
summand of a $\widehat M$--composed space. }

\smallskip
{\it Proof}. Combining Propositions 2(i) and 3(i) with Theorem 4 and Proposition 2, we obtain (i). 
Combining Propositions 3(ii) and 10 with Theorem 4, we obtain the `only if' part of (ii). To prove 
the rest, we observe that every space $\widehat M_n$ is, of course, a complete direct summand of 
the space $\widehat M_\ii$, hence, by Propositions 10 and 11(i), combined with Theorem 3, it is 
metrically projective.  It remains to use Propositions 8 and 11(ii), and then Propositions 3(ii) 
and (again) 10. $\sq$ 

\medskip
{\bf Remark.} Blecher~\cite{ble} considered a different kind of projectivity. This was the operator 
space version of the `lifting property' of some Banach spaces (cf., e.,g.,~\cite[p. 133]{pie}), 
studied in the classical context by Grothendieck~\cite{gro}. This  kind of projectivity  also can 
be treated within the general frame-work of a rigged category and its free objects, but after a 
kind of elaboration of our scheme. Such an approach
%, in a frame-work of the so--called `asymptotic projectivity', 
was used for
%, in the context of 
operator spaces, in~\cite{he6,he7}. As to general matricially normed spaces, this approach leads to 
the following version of Theorem 3.10 in~\cite{ble}: 
 
% In the spirit of~\cite{ble}, we say that a matricially normed space $F$ is almost a direct summand 
% of a matricially normed space $E$, if $F$ is completely isometrically isomorphic to a subspace $G$ of
% $E$, and if for all $\e>0$ there is a completely bounded projection $Q$ of $E$ onto $G$ such that
%  $\|Q\|_{cd}<1+\e$. Then we have:
  
{\it A matricially normed space is projective }(in the just mentioned sense) {\it if, and only if 
it is almost a direct summand of a $\widehat M$--composed space. }

The definition of an almost direct summand repeats word by word the Definition 3.8 in~\cite{ble} 
that was given for operator spaces. 

\bigskip
\centerline{{\bf 7. Some properties of the matricially normed space $\widehat M_n$}}

\bigskip
In this section we again fix a natural $n$.

%For what follows, we shall use the matricially normed spaces $\co_{\min}$ and $\co_{\max}$. The 
%first one is  $\co$ with the matrix-norm, arising after the identifying, for every $m\in{\Bbb N}$, 
%of $M_m(\co)$ with ${\Bbb M}_m$, whereas the second one is $\co$ with the matrix-norm, arising 
%after the identifying of $M_m(\co)$ with ${\Bbb T}_m$ (cf. Section 1). 

\medskip
{\bf Theorem 6}. {\it The underlying normed space of $\widehat M_n$ is ${\Bbb T}_n$}. 

\smallskip
{\it Proof}. Denote the norm on $M_1(\widehat M_n)$ by $\|\cd\|_\bullet$. Take an arbitrary 
element, say $a$, in $M_1(\widehat M_n)$. First, we prove that $\|a\|_\bullet\le\|\cd\|_t$. 
Accordingly, our task is to show that for every  proper couple $(E,v)$ we have 
$\|\va^v(a)\|\le\|a\|_t$. 

Consider
%, in the pure algebraic context, 
the commutative diagram (D2) with $\co$ as $F$ and an arbitrary functional $f$ on $E$ as $\psi$. 
Fix $v\in M_n(E)$ and denote, for brevity, the matrix $f_n(v)\in M_n(\co)$ by $b$. Then for our 
$a$, as for a matrix in $M_n$, we have $f[\va^v(a)]=\va^b(a)$.
%$$
%f[\va^v(a)]=\va^b(a).\eqno(7)
%$$ 
Consequently,  knowing what is $\va^b:M_n\to\co$ (cf. Section 3) and using the latter equation, we 
have $f[\va^v(a)]=tr(ab)$. Therefore the standard duality between ${\Bbb M}_n$ and ${\Bbb T}_n$ 
gives the estimate $|f[\va^v(a)]|\le\|b\|_o\|a\|_t$. 
% $$
%|f[\va^v(a)]|\le\|b\|_o\|a\|_t.\eqno(8)
%$$

Now, using Hahn/Banach theorem, take $f\in E^*;\|f\|=1$ such that $f(\va^v(a))=\|\va^v(a)\|$. Then, 
by the latter inequality, we have 
$$
\|\va^v(a)\|\le\|b\|_o\|a\|_t.\eqno(7)
$$
But it follows from~\cite[Cor. 3.3]{efr} that $f$, being considered as an operator between $E$ and 
$\co_{\min}$, is completely contractive. Since $\|v\|\le1$, this implies that the norm of 
${f_n(v)}$ in $M_n(\co_{\min})={\Bbb M}_n$ is also $\le1$, that is $\|b\|_o\le1$. Therefore the 
needed estimate for $\|\va^v(a)\|$ follows from (7). 
 
Turn to the inverse estimate. By the duality between ${\Bbb M}_n$ and ${\Bbb T}_n$, there exists 
$w\in{\Bbb M}_n;\|w\|_o=1$ such that $tr(aw)=\|a\|_t$. Set $E:=\co_{\min}$ and consider $w$ in the 
unit ball of $M_n({\co_{\min}})$. Then $\va^w:M_n\to\co$ takes $a$ to $\|a\|_t$. Therefore, by (2) 
(with $m=1$), we have $\|a\|_\bullet\ge\|a\|_t$. $\sq$ 

\medskip
Note that the underlying space of the {\it operator} space, playing in the smaller category of 
operator spaces the same role of `building bricks' for free objects, 
%as $\widehat M_n$,%the free object with a one--point base, 
is again ${\Bbb T}_n$~\cite[Prop. 2.7]{he7}. However, our current object, the space $\widehat M_n$, 
is  far away to be an operator space. We have already seen this for $n=1$; now we demonstrate this 
for all $n$. 
%Now we demonstrate this. %More of this, we have the following fact. 

\medskip
Take $p\in[1,\ii)$. A matricially normed space $E$ is said to be {\it $p$--convex} or {\it 
$p$--concave} if for every matrices $u_1,\dots,u_n$ with entries in $E$, we have that 
$\|u_1\oplus...\oplus u_n\|\le\sum_{k=1}^n\|u_k\|^p)^{\frac{1}{p}}$ or $\|u_1\oplus...\oplus 
u_n\|\ge\sum_{k=1}^n\|u_k\|^p)^{\frac{1}{p}}$, respectively. A space that is $p$--convex and 
$p$--concave, is called an $L^p$--space~\cite{efr}. Evidently, an operator space is $p$--convex for 
every $p$. 
%Following ~\cite{efr}, for $p\in[1,\ii)$ we say that a matricially normed space $E$ is an {\it 
%${\cal L}^p$--space}, if for every $E$--valued matrices $u_1,\dots,u_n$ (may be, of different 
%sizes) we have that $\|u_1\oplus...\oplus u_n\|=\sum_{k=1}^n\|u_k\|^p)^{\frac{1}{p}}$. If we 
%replace in the latter expression the sign `=' by $`\le'$ or $`\ge'$, we come to the notion 
%of  a {\it $p$--convex space} or {\it $p$--concave space}, respectively. 

\medskip
{\bf Proposition 12.} {\it The matricially normed space $\widehat M_n$ is not $p$--convex, in 
particular, not an $L^p$--space, for every $p>1$. } 

\smallskip
{\it Proof}. Suppose the contrary. Take any $q$--concave matricially normed space $E\ne0$ with 
$1\le q<p$ (for example, $\co_{\max}$). Then, according to~\cite[Theorem 5.3]{ru2} (cf. 
also~\cite[Prop. 3.3]{colh}), we have $\cc(\widehat M_n,E)=0$. On the other hand, $B_{M_n(E)}$ has 
certainly more than one point. But by virtue of Theorem 2
% because of $\odot_n$--freeness of $\widehat M_n$ , 
there is a bijection between the sets $B_{M_n(E)}$ and $\cc(\widehat M_n,E)$. We came to a 
contradiction. $\sq$ 
%\medskip
%The following proposition shows that the norms of the diagonal matrices in $M_m(\widehat M_n)$ 
%behave similarly to what happens in the case of $L_1$-spaces, but `up to a constant multiplier', 
%namely $\frac{1}{n}$.  

\medskip
{\bf Proposition 13.} {\it Let ${\bf a}\in M_m(\widehat M_n)$ be a block--diagonal matrix ${\bf 
a}=(a_1\oplus...\oplus a_m);a_k\in M_n$. Then $\|{\bf a}\|_m\ge\frac{1}{n}\sum_k\|a_k\|_t$. } 

\smallskip
{\it Proof}. As it is well known, there exist unitary matrices $S_k,T_k: k=1,\dots,m$  such that 
every $S_ka_kT_k$ is a positive diagonal matrix. Note that for ${\bf S}:=(S_1\oplus...\oplus S_m)$ 
and ${\bf T}:=(T_1\oplus...\oplus T_m)$ we have ${\bf S}{\bf a}{\bf T}=(S_1a_1T_1\oplus...\oplus 
S_ma_mT_m)$. Therefore, because of Axiom 2, we can suppose without loss of generality that all  
$a_k$ are positive diagonal matrices. 

Set $E:=\co_{\max}$, so we have $M_m(E)={\Bbb T}_m; m\in{\Bbb N}$. Also set $v:=\frac{1}{n}{\bf 
1}$, where ${\bf 1}$ is the identity matrix in $\widehat M_n$, so we have $v\in B_{M_n(E)}$. It is 
easy to see that $\va^v:M_n\to\co$ takes $a$ to $tr(va)=\frac{1}{n}tr(a)$, hence $(\va^v)_m({\bf 
a})$ is the diagonal $m\times m$--matrix with numbers $\frac{1}{n}tr(a_k)$ on the diagonal. 
Therefore, by (2), we have $\|{\bf a}\|\ge \|(\va^v)_m({\bf 
a})\|_t=\frac{1}{n}\sum_ktr(a_k)=\frac{1}{n}\|a_k\|_t$. $\sq$ 

\medskip
Note that Proposition 12 could be easily deduced from the previous proposition, without applying to 
the triviality of the set $\cc(\widehat M_n,E)$. 

Proposition 13 shows, loosely speaking, that some properties of $\widehat M_n$ resemble (`up to the 
multiplier $\frac{1}{n}$') to those of $L^1$--spaces. Nevertheless we have

\medskip
{\bf Proposition 14.} {\it The matricially normed space $\widehat M_n;n>1$ is not an $L^1$--space. 
} 

\smallskip
{\it Proof}. Suppose the contrary. As a particular case of Proposition 3.2 in~\cite{colh}, every 
functional $f:E\to\co_{\max}$, where $E$ is an $L^1$--space, contractive in the `classical' sense, 
is automatically completely contractive. Since the underlying space of $\widehat M_n$ is ${\Bbb 
T}_n$, this concerns, in particular, $f:\widehat M_n\to\co:a\mapsto tr(a)$. Consequently, the 
operator $f_n:M_n(\widehat M_n)\to M_n(\co_{\max})$ is contractive. In particular, for ${\Bbb I}\in 
M_n(M_n)$ (see Section 5) we have $\|f_n({\Bbb I})\|_t\le\|{\Bbb I}\|$. But, by Proposition 7, we 
have $\|{\Bbb I}\|=1$, and at the same time $f_n({\Bbb I})=\sum_{ij}tr(e_{ij})e_{ji}$ is the 
identity matrix in $M_n$. Therefore $\|f_n({\Bbb I})\|=n>1$, a contradiction. $\sq$

\begin{flushleft}
Moscow State (Lomonosov) University\\ Moscow, 111991, Russia\\
E-mail address: helemskii@rambler.ru 
\end{flushleft}

\end{document}